\documentclass[11pt]{article}
\usepackage{amsthm,amsmath,latexsym,amssymb}

\newcommand{\bP}{{\rm |\kern-.15em P}}
\newcommand{\Q}{\kern.3em\rule{.07em}{.65em}\kern-.3em{\rm Q}}
\newcommand{\R}{{\rm I\kern-.15em R}}
\newcommand{\D}{{\rm |\kern-.15em D}}
\newcommand{\h}{{\rm |\kern-.15em H}}
\newcommand{\C}{\kern.3em\rule{.07em}{.65em}\kern-.3em{\rm C}}
\newcommand{\T}{{\rm T\kern-.35em T}}

\theoremstyle{plain}
\newtheorem{theorem}{Theorem}[section]

\theoremstyle{definition}
\newtheorem{definition}[theorem]{Definition}

\theoremstyle{remark}

\newcommand\blfootnote[1]{%
  \begingroup
  \renewcommand\thefootnote{}\footnote{#1}%
  \addtocounter{footnote}{-1}%
  \endgroup
}

\begin{document}
\title{A disproof of a conjecture of Al Baernstein II}
\author{Ronen Peretz}
 
\maketitle

\begin{abstract}
We answer a question that was asked by Albert Baernstein II, regarding the coefficients of circular symmetrization. The
conjecture is not true generically.
\end{abstract}

\section{The 1974 conjecture of Baernstein}

\blfootnote{\textup{2010} \textit{Mathematics Subject Classification}: \textup{30C55, 30C75, 30C80, 30C85, 30H10}}
\blfootnote{\textit{Key Words and Phrases:} \textup{circular symmetrization, Steiner symmetrization, extremal problems}}

We will recall results from the paper \cite{b}.

\begin{definition}
Let $D$ be a domain in the Riemann sphere $\mathbb{C}\cup\{\infty\}$. The circular symmetrization of $D$ is the 
domain $D^{*}$ that is defined as follows: for each $t\in (0,\infty)$ we define $D(t)=\{\theta\in [0,2\pi]\,|\,te^{i\theta}\in D\}$.
If $D(t)=[0,2\pi]$ then the intersection of $D^{*}$ with the circle $|z|=t$ is the full circle. If $D(t)=\emptyset$
then the intersection of $D^{*}$ with the circle $|z|=t$ is the empty set $\emptyset$. If $D(t)$ is a non trivial subset of
$[0,2\pi]$ which has the measure $|D(t)|=\alpha'$, then the intersection of $D^{*}$ with the circle $|z|=t$ is the unique
circular arc given by $\{te^{i\theta}\,|\,|\theta|< \alpha'/2\}$. Finally $D^{*}$ contain the point $0$ ($\infty$)
if and only if $D$ contains the point $0$ ($\infty$). \\
Let $f\in H(U)$ be one to one and let $F$ be the conformal mapping of $U$ onto $f(U)^{*}$ (both $f(U)$ and $f(U)^{*}$
are simply connected domains) that satisfies $F(0)=|f(0)|$, $F^{'}(0)>0$. Then $F$ is called the circular symmetrization of $f$.
\end{definition}
\noindent
Section (j) of the paper \cite{b} includes a proof of an important principle in symmetrization: \\
Let $f\in H(U)$ and let us denote $D=f(U)$. Let $D_0$ be a simply connected domain that contains $D^{*}$, and let us
assume that $D_0$ is not the full complex plane ($\mathbb{C}$). Let $F$ be a conformal mapping of $U$ onto $D_0$
that satisfies $F(0)=|f(0)|$. The following result is proved in \cite{b}: \\
\\
{\bf Theorem 6. (\cite{b})} {\it
If  $\Phi$ is a convex non-decreasing function on $(-\infty,\infty)$, $f\in H(U)$ and $F$ as above, then for all
$0\le r<1$ we have:
$$
\int_{-\pi}^{\pi}\Phi(\log |f(re^{i\theta})|)d\theta\le\int_{-\pi}^{\pi}\Phi(\log |F(re^{i\theta})|)d\theta.
$$
}
If we choose in Theorem 6 above, $\Phi(x)=e^{2x}$ and assume that we have the following expansions: 
$f(z)=\sum_{n=0}^{\infty}a_{n}z^{n}$ and $F(z)=\sum_{n=0}^{\infty}A_{n}z^{n}$, then we obtain the inequality
$\sum_{n=0}^{\infty}|a_{n}|^{2}r^{2n}\le\sum_{n=0}^{\infty}|A_{n}|^{2}r^{2n}$ for $0\le r<1$. By the definition
of $F$ we have $|A_{0}|=|a_{0}|$, thus if we subtract $|A_{0}|^{2}$ from both sides of the inequality and divide by $r^{2}$
and than take $r\rightarrow 0^{+}$ we obtain $|f'(0)|\le |F'(0)|$, a classical result of Walter Hayman. If $f$
is one-to-one in $U$ then both $D$ and $D^{*}$ are simply connected and we can take $F$ to be a conformal mapping
from $U$ onto $D^{*}$ for which $F(0)=|f(0)|$.

At the end of section (k) in \cite{b} the author asks if the following is true for all $n$: $|a_{n}|\le |A_{n}|$?
Is the following weaker set of inequalities true: $\sum_{k=0}^{n}|a_{k}|^{2}\le\sum_{k=0}^{\infty}|A_{k}|^{2}$?
A. Baernstein II, remarks that these last inequalities if true, would prove a conjecture of Littlewood: If $f$
is one-to-one and analytic in $U$ and if $f(z)\ne 0$, for $z\in U$, then for each $n>1$ we have: $a_{n}\le 4n|a_{0}|$. \\
It was proved in \cite{ds} that the Bieberbach's conjecture implies the above Littlewood's conjecture. Since by now we know
the Bieberbach's conjecture to be true, \cite{d}, Littlewood's conjecture is true as well.

\section{A disproof of the conjecture of Al Baernstein II}
Concerning the first question posed by Albert Baernstein II (above) we prove the following:

\begin{theorem}\label{thm1}
If $f(z)=\sum_{n=0}^{\infty}a_{n}z^{n}$ is analytic, one-to-one in $U$ and if 
$F(z)=\sum_{n=0}^{\infty}A_{n}z^{n}$ is the circular symmetrization of $f(z)$, then we have: 
$\sum_{n=0}^{\infty}n|a_{n}|^{2}=\sum_{n=0}^{\infty}n|A_{n}|^{2}$ and either for all $n=0,1,2,\ldots $
we have $|a_{n}|=|A_{n}|$, or there exist $1\le n_{1},n_{2}$ such that $|a_{n_{1}}|<|A_{n_{1}}|$ and
$|A_{n_{2}}|<|a_{n_{2}}|$.
\end{theorem} 
\noindent
Theorem \ref{thm1} answers the problem mentioned above that was raised by Albert Baernstein II. The answer in negative. \\
\\
{\bf Proof.} \\
Let $f$ be a conformal function defined on $U$. Here conformal means analytic and one-to-one.
We assume that $f(U)$ has a finite area (otherwise we replace $f(z)$ by $f(rz)$ for $0<r<1$). Let us denote $D=f(U)$, 
and let $F$ be a conformal mapping of $U$ onto the symmetrization $D^{*}$ such that $F(0)=|f(0)|$. Let us denote by $S(D)$ 
and by $S(D^{*})$ the areas of the respective domains. We will use $td\phi\cdot dt$ for the area element in polar coordinates. 
Then we have the identities:
$$
S(D)=\int_{0}^{\infty}\int_{D(t)}td\phi\cdot dt=\int_{0}^{\infty}t|D(t)|dt,\,\,\,\,S(D^{*})=\int_{0}^{\infty}t|D^{*}(t)|dt.
$$
By the definition of $D^{*}$ it follows that for all $0\le t<\infty$ we have $D(t)=D^{*}(t)$ and hence $S(D)=S(D^{*})$
(the well-known fact that circular symmetrization is an area preserving transformation). On the other hand we have
$S(D)=\int_{0}^{1}\int_{0}^{2\pi}r|f'(re^{i\theta})|^{2}d\theta dr$ and if $f(z)=\sum_{n=0}^{\infty}a_{n}z^{n}$ and
$F(z)=\sum_{n=0}^{\infty}A_{n}z^{n}$ then we obtain the well known formulas: $S(D)=\pi\sum_{n=0}^{\infty}n|a_{n}|^{2}$, and
$S(D^{*})=\pi\sum_{n=0}^{\infty}n|A_{n}|^{2}$. We conclude that $\sum_{n=0}^{\infty}n|a_{n}|^{2}=\sum_{n=0}^{\infty}
n|A_{n}|^{2}$. We recall that by the definition of $F$ we have $A_{0}=|a_{0}|$ and by Hayman's result (see \cite{ha}) $|a_{1}|\le |A_{1}|$
and so either $|a_{n}|=|A_{n}|$ for $n=0,1,2,\ldots $ or there exist $1\le n_{1},n_{2}$ so that $|a_{n_{1}}|<|A_{n_{1}}|$,
and $|A_{n_{2}}|<|a_{n_{2}}|$. In the first case (in which we have equalities of the absolute values for all the
Taylor coefficients) the uniqueness result in \cite{ha} completes our proof. Hayman refers to a result of Jenkins,
\cite{j1} and to his own paper \cite{ha1}. Using the results in those papers one can show that strict inequality holds
in $|a_{1}|<|A_{1}|$ unless $f(U)=f(U)^{*}$ and $f(z)=F(e^{i\lambda}z)$ for some real $\lambda$.
It is worth mentioning here related uniqueness results of Jenkins (in \cite{j2}) and of Ess\'en and Shea
(in \cite{es}). Our proof is now complete. $\qed $

\noindent
{\it Ronen Peretz \\
Department of Mathematics \\ Ben Gurion University of the Negev \\
Beer-Sheva , 84105 \\ Israel \\ E-mail: ronenp@math.bgu.ac.il} \\ 
 
\end{document}